\documentclass{article}
\usepackage[T1]{fontenc}
\usepackage{lmodern}
\usepackage{amsmath,amssymb,amsthm,mathtools}
\usepackage{fullpage}
\usepackage{microtype}
\usepackage{tikz}
\usetikzlibrary{arrows.meta}
\usepackage[hidelinks]{hyperref}
\hypersetup{pdftitle={Tableau shifts for graded multiplicities in the (GL4,O4) Kostant-Rallis setting},pdfauthor={Andrew Frohmader and Jeb Willenbring}}

\newtheorem{thm}{Theorem}[section]
\newtheorem{lemma}[thm]{Lemma}
\newtheorem{cor}[thm]{Corollary}
\newtheorem{prop}[thm]{Proposition}
\theoremstyle{definition}
\newtheorem{defn}[thm]{Definition}

\DeclareMathOperator{\GL}{GL}
\DeclareMathOperator{\SL}{SL}
\DeclareMathOperator{\SOr}{SO}
\DeclareMathOperator{\Or}{O}

\newcommand{\Ha}{\mathcal H}

\title{Tableau shifts for graded multiplicities in the $(\GL_4,\Or_4)$ Kostant-Rallis setting}
\author{Andrew Frohmader \and Jeb Willenbring}
\date{November 2026}

\begin{document}
\maketitle

\begin{abstract}
Wallach and Willenbring computed the graded multiplicities of irreducible $\SOr_4$-representations in the coordinate ring of the nilpotent symmetric $4\times4$ matrices. Using a tableau formula for the symmetric pair $(\GL_4,\Or_4)$, we construct two explicit injections that explain shift formulas for these multiplicities allowing us to give a tableau proof of the Wallach-Willenbring shift and a companion recurrence in a second direction on the dominant weight lattice of $\SOr_4$. Iterating the recurrences gives a positive expression for the graded multiplicities.
\end{abstract}

\section{Introduction}
Let $\mathcal N_\theta$ be the variety of nilpotent symmetric $4\times4$ complex matrices. We study the graded multiplicities of irreducible $\Or_4$- and $\SOr_4$-representations in $\mathbb C[\mathcal N_\theta]$, with the grading induced by polynomial degree and the action induced by conjugation. By the Kostant-Rallis theorem \cite{KostantRallis1971}, this coordinate ring is isomorphic, as a graded representation, to the space $\Ha$ of harmonic polynomials for the symmetric pair $(\GL_4,\Or_4)$.

This example also arises in the decomposition of curvature tensors in dimension four. As an $\Or_4$-representation, the $20$-dimensional space of algebraic curvature tensors on $\mathbb C^4$ decomposes into the scalar, traceless Ricci, and Weyl curvature components, of dimensions $1$, $9$, and $10$, respectively \cite[\S10.3]{GoodmanWallach2009}. The scalar component is trivial, while the pullback of the Weyl curvature component along the double cover $\SL_2\times\SL_2\to\SOr_4$ is the direct sum of two binary quartic representations, whose symmetric powers are classically understood. The traceless Ricci component is isomorphic to the space of traceless symmetric $4\times4$ matrices under conjugation, so $\mathcal N_\theta$ is the null cone in this component. 

Wallach and Willenbring determined all of the $\SOr_4$ graded multiplicities, giving a shift recurrence and complete parity-dependent formulas \cite[Theorem~7]{WallachWillenbring2000}. Their verification used MAPLE; Willenbring also gave a proof of the shift using a virtual-character identity, spherical harmonics, and Clebsch-Gordan \cite[\S6.6]{Willenbring2000}. Our purpose is to give a direct tableau explanation of the recurrences and their $\Or_4$ refinements, using \cite[Theorem~6.1]{Frohmader2026}. The resulting companion recurrence and positive closed form reorganize the known multiplicities; the explicit injections and their complements explain the positive correction terms.

For $n\geq0$, write
\[
[n]_q=\frac{1-q^n}{1-q},\qquad [n]_q!=\prod_{j=1}^{n}[j]_q,\qquad
\begin{bmatrix}n\\m\end{bmatrix}_q=\frac{[n]_q!}{[m]_q![n-m]_q!}\quad(0\leq m\leq n).
\]
Thus $[0]_q=0$ and $[0]_q!=1$. We set $\begin{bmatrix}n\\m\end{bmatrix}_q=0$ when $m<0$ or $m>n$, and all empty sums are zero.

Let $\varpi_1=\frac12(\epsilon_1-\epsilon_2)$ and $\varpi_2=\frac12(\epsilon_1+\epsilon_2)$ be the fundamental weights for $\mathfrak{so}_4$. Since the central element $-I$ acts trivially by conjugation, only representations descending to $\SOr_4/\{\pm I\}$ can occur. Their highest weights have the form $2k\varpi_1+2\ell\varpi_2$, with $k,\ell\geq0$. We write $p_{k,\ell}$ for the graded multiplicity of the corresponding irreducible representation in $\mathbb C[\mathcal N_\theta]$.

\subsection{Two shift formulas}
Restriction from $\Or_4$, recalled in Section~\ref{subsec:O_4-model}, gives $p_{k,\ell}=p_{\ell,k}$. It therefore suffices to work on the half-lattice $k\geq\ell\geq0$. The first recurrence is the Wallach-Willenbring shift \cite[p.~445]{WallachWillenbring2000}.

\begin{thm}\label{thm:intro-WW-shift}
For $k\geq\ell\geq0$,
\[
p_{k+2,\ell}-q^2p_{k,\ell}=q^{2k-\ell+4}[2\ell+1]_q.
\]
\end{thm}

Our tableau proof appears at the end of Section~\ref{sec:2-2-shift}. The second recurrence increments both weight coordinates.

\begin{thm}\label{thm:intro-2-shift}
For $k\geq\ell\geq0$,
\[
p_{k+1,\ell+1}-qp_{k,\ell}=q^{k+\ell+2}\bigl([k+2]_q+q^2[\ell]_q\bigr).
\]
\end{thm}

\begin{figure}[htbp]
\centering
\begin{tikzpicture}[
    x=0.62cm,
    y=0.62cm,
    every node/.style={font=\scriptsize},
    shiftarrow/.style={
        -{Stealth[length=1.8mm,width=1.3mm]},
        draw=black!75,
        line width=0.75pt,
        shorten <=2.2pt,
        shorten >=2.2pt
    },
    evenpt/.style={fill=black,draw=black},
    oddpt/.style={fill=white,draw=black,line width=0.45pt}
]
    \def\N{5}

    \draw[->,line width=0.6pt] (-0.18,0) -- (\N+0.42,0) node[right] {$k$};
    \draw[->,line width=0.6pt] (0,-0.18) -- (0,\N+0.42) node[above] {$\ell$};

    \foreach \x in {0,...,\N} {
        \draw (\x,0.055) -- (\x,-0.055);
        \node[below=2pt] at (\x,0) {\x};
    }
    \foreach \y in {1,...,\N} {
        \draw (0.055,\y) -- (-0.055,\y);
        \node[left=2pt] at (0,\y) {\y};
    }

    \draw[densely dashed,black!35,line width=0.55pt] (0,0) -- (\N,\N);

    \foreach \k in {0,...,\N} {
        \foreach \l in {0,...,\k} {
            \pgfmathtruncatemacro{\parity}{mod(\k-\l,2)}
            \ifnum\parity=0
                \filldraw[evenpt] (\k,\l) circle (1.7pt);
            \else
                \filldraw[oddpt] (\k,\l) circle (1.7pt);
            \fi
        }
    }

    \draw[shiftarrow] (2,1) -- (3,2);
    \draw[shiftarrow] (4,1) -- (5,2);

    \draw[shiftarrow] (2,1) .. controls (2.65,1.34) and (3.35,1.34) .. (4,1);
    \draw[shiftarrow] (3,2) .. controls (3.65,2.34) and (4.35,2.34) .. (5,2);

    \node[fill=white,inner sep=1.2pt,anchor=east] at (2.36,1.67) {$A$};
    \node[fill=white,inner sep=1.2pt,anchor=north] at (3,0.82) {$B$};

    \begin{scope}[shift={(6.1,3.4)}]
        \filldraw[evenpt] (0,0) circle (1.7pt);
        \node[anchor=west] at (0.2,0) {$k-\ell\equiv0\pmod2$};

        \filldraw[oddpt] (0,-0.7) circle (1.7pt);
        \node[anchor=west] at (0.2,-0.7) {$k-\ell\equiv1\pmod2$};

        \node[anchor=west] at (0,-1.65) {$A=(1,1)$};
        \node[anchor=west] at (0,-2.25) {$B=(2,0)$};
    \end{scope}
\end{tikzpicture}
\caption{The half-lattice $k\geq\ell\geq0$ and a commuting square for the translations $A$ and $B$.}
\label{fig:SO4-half-lattice}
\end{figure}
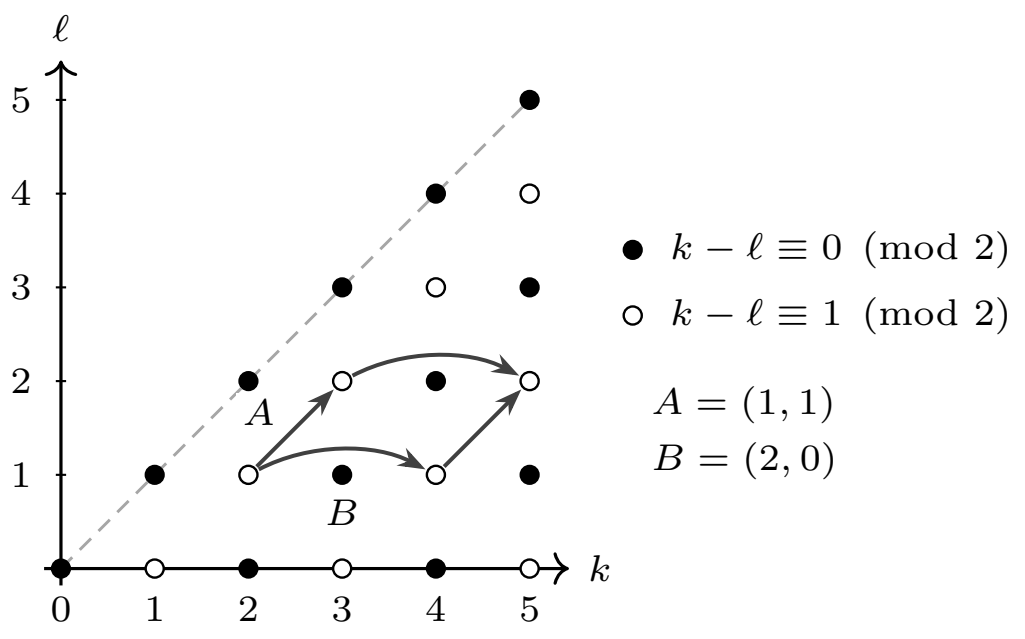

This is proved at the end of Section~\ref{sec:2-shift}. The two recurrences correspond to the translations
\[
A(k,\ell)=(k+1,\ell+1),\qquad B(k,\ell)=(k+2,\ell).
\]
They preserve the parity of $k-\ell$ and generate the two parity classes from $(0,0)$ and $(1,0)$, respectively. Figure~\ref{fig:SO4-half-lattice} shows a commuting square. In Section~\ref{sec:closed-forms}, we iterate the recurrences to obtain a positive closed form.

\subsection{Restricting from \texorpdfstring{$\Or_4$ to $\SOr_4$}{O4 to SO4}}\label{subsec:O_4-model}
We index polynomial irreducible $\Or_4$-representations by
\[
\widehat{\Or}_4=\{\lambda:\text{the first two columns of $\lambda$ contain at most four boxes}\},
\]
where $\lambda$ is a Young diagram. Write $m_\lambda$ for the graded multiplicity of the representation indexed by $\lambda$ in $\mathbb C[\mathcal N_\theta]$, and abbreviate $m_{(a,b,c,d)}$ to $m_{a,b,c,d}$, omitting trailing zeros.

The restriction rule \cite[Theorem~5.5.24]{GoodmanWallach2009} gives
\[
\begin{aligned}
p_{k,\ell}&=m_{k+\ell,k-\ell}=p_{\ell,k} && (k>\ell\geq0),\\
p_{k,k}&=m_{2k}+m_{2k,1,1} && (k\geq1),\\
p_{0,0}&=m_0.
\end{aligned}
\]
Indeed, $(a,b)$ with $b>0$ restricts to the two highest weights $a\epsilon_1\pm b\epsilon_2$, while $(a)$ and $(a,1,1)$ have the same irreducible restriction, of highest weight $a\epsilon_1$. The determinant representation $(1,1,1,1)$ has zero multiplicity here.

Away from the diagonal, the translations $A$ and $B$ are described at the $\Or_4$ level by the Young diagram shifts $(a,b)\mapsto(a+2,b)$ and $(a,b)\mapsto(a+2,b+2)$, respectively. The diagonal cases are treated separately.

\section{Tableau model and shift maps}\label{sec:tableaux}
We read a tableau row by row from bottom to top, reading each row from left to right. For $1\leq i\leq3$, the \emph{signature rule} computes $\varphi_i$ from this word: put a $-$ over each $i$ and a $+$ over each $i+1$, then successively cancel adjacent $(+,-)$ pairs, ignoring blanks. The number of remaining $-$'s is $\varphi_i$.

\begin{defn}
A semistandard tableau $T$ with entries in $\{1,2,3,4\}$ is an \emph{$\Or_4$-tableau} if, for $1\leq i\leq4$,
\[
\#\{\text{boxes in the first two columns of $T$ with entry $\leq i$}\}\leq i.
\]
Let $\mathcal T^\lambda$ denote the set of these tableaux of shape $\lambda$ in which each entry occurs an even number of times. We call this even-content requirement the \emph{zero-$M$-weight} condition, see \cite{Frohmader2026}.
\end{defn}

No $1$ occurs in $T\in\mathcal T^\lambda$: column strictness confines all $1$'s to the top row, and the even-content condition would force its first two entries to be $1$, violating the $\Or_4$ condition. The second row therefore contains only $3$'s and $4$'s, with at most one $3$, since two second-row $3$'s would put four entries at most $3$ in the first two columns. A third row must consist of a single $4$, and then the second row is a single $3$; any further box in these rows would violate the condition on the first two columns. Hence every tableau has one of the forms
\[
\begin{array}{c@{\qquad\text{or}\qquad}c}
\begin{array}{l}2^A\,3^B\,4^C\\3^D\,4^E\end{array}
&\begin{array}{l}2^A\,3^B\,4^C\\3\\4\end{array}
\end{array}
\qquad\text{with }D\leq1.
\]
The first form includes one-row shapes. Since $\varphi_1(T)=0$, the formula in \cite[Theorem~6.1]{Frohmader2026} specializes as follows.

\begin{thm}\label{thm:graded-mults}
For $\lambda\in\widehat{\Or}_4$,
\[
m_\lambda=\sum_{T\in\mathcal T^\lambda}q^{d(T)},\qquad
d(T)=\frac{|\lambda|}{2}+\lceil\varphi_2(T)\rceil_{\mathrm{ev}}+\frac{\lceil\varphi_3(T)\rceil_{\mathrm{ev}}}{2},
\]
where $\lceil x\rceil_{\mathrm{ev}}$ is the smallest even integer greater than or equal to $x$.
\end{thm}

For $a \ge b \ge c \ge 0$, define the \emph{$(2)$-shift}
\[
S_1:\mathcal T^{(a,b,c)}\hookrightarrow\mathcal T^{(a+2,b,c)}
\]
by appending two $4$'s to the first row. Define the \emph{$(2,2)$-shift}
\[
S_2:\mathcal T^{(a,b)}\hookrightarrow\mathcal T^{(a+2,b+2)},\qquad
S_2\left(\begin{array}{l}2^A\,3^B\,4^C\\3^D\,4^E\end{array}\right)
=\begin{array}{l}2^A\,3^{B+2}\,4^C\\3^D\,4^{E+2}.\end{array}
\]
It is easy to see that both maps preserve the $\Or_4$ condition, preserve the zero-$M$-weight condition, and are injective. A tableau lies in the image of $S_1$ exactly when its top row contains at least two $4$'s. A tableau lies in the image of $S_2$ exactly when its top row contains at least two $3$'s and its bottom row contains at least two $4$'s.

\begin{prop}\label{prop:tableau-shifts}
The maps $S_1$ and $S_2$ satisfy
\[
d(S_1(T))=d(T)+1,\qquad d(S_2(T))=d(T)+2.
\]
Furthermore,
\[
m_{a+2,b,c}-qm_{a,b,c}\in\mathbb Z_{\geq0}[q],\qquad
m_{a+2,b+2}-q^2m_{a,b}\in\mathbb Z_{\geq0}[q],
\]
and each difference is the sum of $q^{d(T)}$ over the tableaux outside the corresponding image.
\end{prop}

\begin{proof}
We first prove the statements about $d$. For $S_1$, the added $4$'s can only affect the $3$-signature. They append two $+$'s that cancel no $-$'s, so $\varphi_3$ is unchanged. For a two-row tableau in the displayed form, the signature rule gives
\[
\varphi_1(T)=0,\qquad \varphi_2(T)=\max(0,A-D),\qquad \varphi_3(T)=D+\max(0,B-E).
\]
These expressions are preserved by $S_2$, which increases both $B$ and $E$ by two. The degree changes follow from the increases in size. The multiplicity assertions follow from Theorem~\ref{thm:graded-mults} and the injectivity of $S_1$ and $S_2$.
\end{proof}

\section{The \texorpdfstring{$(2)$}{(2)}-shift}\label{sec:2-shift}
The $(2)$-shift corresponds to the translation $A$. We treat one-row, hook, and two-row shapes separately.

\subsection{One-row shapes}
\begin{thm}\label{thm:one-row-(2)-recurrence}
For even $a\geq0$,
\[
m_{a+2}-qm_a=q^{a+2}\left[\frac a2+2\right]_q.
\]
\end{thm}

\begin{proof}
Any $T \in \mathcal{T}^{(a+2)} \setminus S_1(\mathcal{T}^{(a)})$ has at most one $4$, hence none by the zero-$M$-weight condition. The complementary tableaux are therefore
\[
T_i=2^{2i}\,3^{a+2-2i},\qquad 0\leq i\leq\frac a2+1.
\]
The signature rule gives $\lceil\varphi_2(T_i)\rceil_{\mathrm{ev}}=2i$ and $\lceil\varphi_3(T_i)\rceil_{\mathrm{ev}}=a+2-2i$, so $d(T_i)=a+2+i$. Thus
\[
m_{a+2}-qm_a=\sum_{i=0}^{a/2+1}q^{a+2+i}=q^{a+2}\left[\frac a2+2\right]_q.
\]
\end{proof}

\begin{cor}\label{cor:one-row-closed-form}
For $n\geq0$,
\[
m_{2n}=q^n\begin{bmatrix}n+2\\2\end{bmatrix}_q.
\]
\end{cor}

\begin{proof}
Starting from $m_0=1$, apply Theorem~\ref{thm:one-row-(2)-recurrence} and the $q$-Pascal identity
\[
\begin{bmatrix}n+2\\2\end{bmatrix}_q=\begin{bmatrix}n+1\\2\end{bmatrix}_q+q^n[n+1]_q
\]
to obtain the result by induction.
\end{proof}

\subsection{Hook shapes}
For even $a\geq2$, every $T\in\mathcal T^{(a+2,1,1)}$ has the form
\[
T=\begin{array}{l}2^{2A+2}\,3^{2B+1}\,4^{2C+1}\\3\\4\end{array},\qquad
A,B,C\geq0.
\]
Deleting the two lower boxes and two $2$'s, one $3$, and one $4$ from the top row defines a bijection
\[
\Psi:\mathcal T^{(a+2,1,1)}\longrightarrow\mathcal T^{(a-2)},\qquad
\Psi(T)=2^{2A}\,3^{2B}\,4^{2C}.
\]

\begin{thm}\label{thm:hook-(2)-recurrence}
For even $a\geq2$,
\[
m_{a+2,1,1}=q^6m_{a-2},\qquad
m_{a+2,1,1}-qm_{a,1,1}=q^{a+4}\left[\frac a2\right]_q.
\]
\end{thm}

\begin{proof}
The rounded signature values of $T$ are $\lceil\varphi_2(T)\rceil_{\mathrm{ev}}=2A+2$ and $\lceil\varphi_3(T)\rceil_{\mathrm{ev}}=2B+2$, while those of $\Psi(T)$ are $2A$ and $2B$. In addition, the change in size contributes 3. So, 
\[
d(T)-d(\Psi(T))=3+2+1=6,
\]
proving the first identity. For the recurrence, $a=2$ follows from $m_{2,1,1}=0$ and $m_{4,1,1}=q^6$; for $a\geq4$, apply Theorem~\ref{thm:one-row-(2)-recurrence} and $m_{a+2,1,1}=q^6m_{a-2}$.
\end{proof}

\begin{cor}\label{cor:hook-closed-form}
For $n\geq1$,
\[
m_{2n,1,1}=q^{n+4}\begin{bmatrix}n\\2\end{bmatrix}_q.
\]
\end{cor}

\begin{proof}
For $n=1$, both sides vanish. For $n\geq2$, combine $m_{2n,1,1}=q^6m_{2n-4}$ with Corollary~\ref{cor:one-row-closed-form}.
\end{proof}

\subsection{Two-row shapes}
We use the auxiliary family
\[
\mathcal F_{2u+\epsilon,2v+\epsilon}\coloneqq
\left\{T\in\mathcal T^{(2u+\epsilon,2v+\epsilon)}:
\begin{array}{l}\text{the bottom row consists entirely of $4$'s, and}\\\text{the top row contains at most one $4$}\end{array}\right\},
\]
where $u\geq v\geq0$, $\epsilon\in\{0,1\}$, and $u-v\geq\epsilon$. This last condition excludes equal odd row lengths where column strictness prevents a top-row $4$, leaving an odd total number of $4$'s. Set
\[
F_{2u+\epsilon,2v+\epsilon}=\sum_{T\in\mathcal F_{2u+\epsilon,2v+\epsilon}}q^{d(T)}.
\]

\begin{lemma}\label{lemma:2-row-(2)-shift-lemma}
Let $u\geq v\geq0$, $\epsilon\in\{0,1\}$, and $u-v\geq\epsilon$. Then
\[
F_{2u+\epsilon,2v+\epsilon}=q^{2u+\epsilon}[u-v+1]_q+q^{3u-v+\epsilon+2}[v]_{q^2}.
\]
\end{lemma}

\begin{proof}
The tableaux in this family are
\[
T_i=\begin{array}{l}2^{2i}\,3^{2u-2i}\,4^\epsilon\\4^{2v+\epsilon}\end{array},\qquad 0\leq i\leq u.
\]
Here $\lceil\varphi_2(T_i)\rceil_{\mathrm{ev}}=2i$ and $\lceil\varphi_3(T_i)\rceil_{\mathrm{ev}}=2\max(0,u-v-i)$, so
\[
d(T_i)=\begin{cases}
2u+\epsilon+i,&0\leq i\leq u-v,\\
u+v+\epsilon+2i,&u-v<i\leq u.
\end{cases}
\]
The two ranges contribute $q^{2u+\epsilon}[u-v+1]_q$ and $q^{3u-v+\epsilon+2}[v]_{q^2}$, respectively.
\end{proof}

\begin{thm}\label{thm:two-row-(2)-recurrence}
Let $a\geq b>0$ and $a\equiv b\pmod2$. Then
\[
m_{a+2,b}-qm_{a,b}=q^{a+2}\left(\left[\frac{a+b}{2}+2\right]_q+q^2\left[\frac{a-b}{2}\right]_q\right).
\]
\end{thm}

\begin{proof}
Write $a=2u+\epsilon$, $b=2v+\epsilon$, with $\epsilon\in\{0,1\}$. The tableaux in $\mathcal T^{(a+2,b)}\setminus S_1(\mathcal T^{(a,b)})$ have at most one top-row $4$. Those with bottom row consisting entirely of $4$'s contribute $F_{2u+\epsilon+2,2v+\epsilon}$. Every remaining tableau has the form
\[
T=\begin{array}{l}2^{A+2}\,3^{B+1}\,4^{1-\epsilon}\\3\,4^{C+1-\epsilon}\end{array},\qquad
A+B=2u-2+2\epsilon,\quad C=2v-2+2\epsilon,
\]
where $A,B,C$ are nonnegative even integers. The map
\[
\Phi(T)=\begin{array}{l}2^A\,3^B\\4^C\end{array}
\]
is a bijection from these tableaux to $\mathcal F_{2u-2+2\epsilon,2v-2+2\epsilon}$. Its target parameters are nonnegative: $b>0$ gives $v\geq1$ if $\epsilon=0$, and $v\geq0$ if $\epsilon=1$. The signature rule gives
\[
\begin{aligned}
\lceil\varphi_2(T)\rceil_{\mathrm{ev}}&=A+2,&\lceil\varphi_2(\Phi(T))\rceil_{\mathrm{ev}}&=A,\\
\lceil\varphi_3(T)\rceil_{\mathrm{ev}}&=2+\max(0,B-C),&\lceil\varphi_3(\Phi(T))\rceil_{\mathrm{ev}}&=\max(0,B-C).
\end{aligned}
\]
Since $|T|-|\Phi(T)|=6-2\epsilon$, we obtain $d(T)-d(\Phi(T))=6-\epsilon$. Consequently,
\[
m_{a+2,b}-qm_{a,b}=F_{2u+\epsilon+2,2v+\epsilon}+q^{6-\epsilon}F_{2u-2+2\epsilon,2v-2+2\epsilon}.
\]
By Lemma~\ref{lemma:2-row-(2)-shift-lemma}, this equals
\[
q^{2u+\epsilon+2}\left([u-v+2]_q+q^2[u-v+1]_q+q^{u-v+3}\bigl([v]_{q^2}+q[v+\epsilon-1]_{q^2}\bigr)\right).
\]
The identities
\[
\begin{aligned}
[v]_{q^2}+q[v+\epsilon-1]_{q^2}&=[2v+\epsilon-1]_q,\\
[u-v+2]_q+q^2[u-v+1]_q&=[u-v+3]_q+q^2[u-v]_q,\\
[u-v+3]_q+q^{u-v+3}[2v+\epsilon-1]_q&=[u+v+\epsilon+2]_q
\end{aligned}
\]
give the result.
\end{proof}

\begin{proof}[Proof of Theorem~\ref{thm:intro-2-shift}]
We use the restriction rules in Section \ref{subsec:O_4-model}. At $k=\ell=0$, Theorem~\ref{thm:one-row-(2)-recurrence} and $m_{2,1,1}=0$ give $p_{1,1}-qp_{0,0}=q^2[2]_q$. For $k\geq1$, Theorem~\ref{thm:one-row-(2)-recurrence}, and Theorem~\ref{thm:hook-(2)-recurrence} give
\[
\begin{aligned}
p_{k+1,k+1}-qp_{k,k}
&=(m_{2k+2}-qm_{2k})+(m_{2k+2,1,1}-qm_{2k,1,1})\\
&=q^{2k+2}[k+2]_q+q^{2k+4}[k]_q.
\end{aligned}
\]
For $k>\ell\geq0$, Theorem~\ref{thm:two-row-(2)-recurrence} yields
\[
\begin{aligned}
p_{k+1,\ell+1}-qp_{k,\ell}
&=m_{k+\ell+2,k-\ell}-qm_{k+\ell,k-\ell}\\
&=q^{k+\ell+2}\bigl([k+2]_q+q^2[\ell]_q\bigr).
\end{aligned}
\]
\end{proof}

\section{The \texorpdfstring{$(2,2)$}{(2,2)}-shift}\label{sec:2-2-shift}
The $(2,2)$-shift corresponds to the translation $B$. We first treat genuine two-row source shapes and then the one-row boundary.

\begin{thm}\label{thm:two-row-(2,2)-shift}
Let $a\geq b>0$ with $a\equiv b\pmod2$. Then
\[
m_{a+2,b+2}-q^2m_{a,b}=q^{\frac{a+3b}{2}+4}[a-b+1]_q.
\]
\end{thm}

\begin{proof}
A tableau lies in the image of $S_2$ exactly when its top row contains at least two $3$'s and its bottom row contains at least two $4$'s. Since $b>0$, the target bottom row has length at least three and at most one $3$, so the second condition is automatic. Outside the image, the number of top-row $3$'s is therefore $\delta\in\{0,1\}$; the zero-$M$-weight condition forces the same number in the bottom row.

Writing $i$ for the number of top-row entries less than $4$, the complementary tableaux are precisely
\[
T_i=\begin{array}{l}2^{i-\delta}\,3^\delta\,4^{a+2-i}\\3^\delta\,4^{b+2-\delta}\end{array},\qquad
b+2\leq i\leq a+2,\quad\delta\equiv i\pmod2.
\]
The lower bound is column strictness, and the parity condition ensures zero-$M$-weight. The signature rule gives
\[
\lceil\varphi_2(T_i)\rceil_{\mathrm{ev}}=i-\delta,\qquad
\lceil\varphi_3(T_i)\rceil_{\mathrm{ev}}=2\delta,
\]
so $d(T_i)=\frac{a+b}{2}+2+i$. Thus
\[
m_{a+2,b+2}-q^2m_{a,b}=\sum_{i=b+2}^{a+2}q^{\frac{a+b}{2}+2+i}=q^{\frac{a+3b}{2}+4}[a-b+1]_q.
\]
\end{proof}

At the one-row boundary, we have

\begin{thm}\label{thm:one-row-(2,2)-shift}
Let $a\geq0$ be even. Then
\[
m_{a+2,2}-q^2m_a=q^{\frac a2+4}\left[\frac a2+1\right]_{q^2}+q^{\frac a2+5}\begin{bmatrix}\frac a2+1\\2\end{bmatrix}_q.
\]
\end{thm}

\begin{proof}
The bottom row of a tableau in $\mathcal T^{(a+2,2)}$ is either $34$ or $44$. Every tableau with bottom row $34$ is outside the image of $S_2$. Moving its bottom-row $4$ to a third row gives a bijection with $\mathcal T^{(a+2,1,1)}$. Both tableaux have top row
\[
2^{2r+2}\,3^{2s+1}\,4^{2t+1},\qquad r,s,t\geq0,
\]
and the same signature statistics $\varphi_2=2r+1$, $\varphi_3=2s+1$. Thus the bijection preserves degree, and these tableaux contribute $m_{a+2,1,1}$. When $a=0$, both sets are empty.

With bottom row $44$, the tableaux outside the image of $S_2$ have no $3$'s and are precisely
\[
T_i=\begin{array}{l}2^{2+2i}\,4^{a-2i}\\44\end{array},\qquad 0\leq i\leq\frac a2.
\]
Here $\varphi_2(T_i)=2+2i$ and $\varphi_3(T_i)=0$, so $d(T_i)=\frac a2+4+2i$. Consequently,
\[
m_{a+2,2}-q^2m_a=m_{a+2,1,1}+q^{\frac a2+4}\left[\frac a2+1\right]_{q^2}.
\]
Corollary~\ref{cor:hook-closed-form} gives the stated expression.
\end{proof}

\begin{proof}[Proof of Theorem~\ref{thm:intro-WW-shift}]
For $k>\ell\geq0$, Section~\ref{subsec:O_4-model} and Theorem~\ref{thm:two-row-(2,2)-shift} yield
\[
\begin{aligned}
p_{k+2,\ell}-q^2p_{k,\ell}
&=m_{k+\ell+2,k-\ell+2}-q^2m_{k+\ell,k-\ell}\\
&=q^{2k-\ell+4}[2\ell+1]_q.
\end{aligned}
\]
On the diagonal, $k=0$ follows from Theorem~\ref{thm:one-row-(2,2)-shift}, which gives $p_{2,0}-q^2p_{0,0}=q^4$. For $k\geq1$, Theorem~\ref{thm:one-row-(2,2)-shift} and Corollary~\ref{cor:hook-closed-form} give
\[
\begin{aligned}
p_{k+2,k}-q^2p_{k,k}
&=m_{2k+2,2}-q^2(m_{2k}+m_{2k,1,1})\\
&=q^{k+4}[k+1]_{q^2}+q^{k+5}\begin{bmatrix}k+1\\2\end{bmatrix}_q-q^{k+6}\begin{bmatrix}k\\2\end{bmatrix}_q\\
&=q^{k+4}[2k+1]_q.
\end{aligned}
\]
The last equality follows from
\[
\begin{bmatrix}k+1\\2\end{bmatrix}_q-q\begin{bmatrix}k\\2\end{bmatrix}_q=[k]_{q^2},\qquad
[k+1]_{q^2}+q[k]_{q^2}=[2k+1]_q.
\]
\end{proof}

\section{Closed forms}\label{sec:closed-forms}
The two recurrences determine all $p_{k,\ell}$ from
\[
p_{0,0}=1,\qquad p_{1,0}=0.
\]
Indeed, the empty tableau gives $p_{0,0}=m_0=1$, represented by the constant functions, while $p_{1,0}=m_{1,1}=0$ because $\mathcal T^{(1,1)}$ is empty. We first iterate in the $B$-direction along $\ell=0$, then in the $A$-direction, see Figure~\ref{fig:SO4-half-lattice}. Wallach-Willenbring instead begin with formulas for $k=\ell$ and $k=\ell+1$ and apply the $B$-shift.

\begin{lemma}\label{lemma:boundary-f-values}
For $t\geq0$,
\[
p_{2t,0}=q^{2t}[t+1]_{q^2},\qquad p_{2t+1,0}=q^{2t+4}[t]_{q^2}.
\]
\end{lemma}

\begin{proof}
We prove both formulas simultaneously by induction on \(t\). For \(t=0\),
\[
p_{0,0}=1=[1]_{q^2},
\qquad
p_{1,0}=0=q^4[0]_{q^2}.
\]
Assume the formulas hold for \(t\). By Theorem~\ref{thm:intro-WW-shift},
\[
\begin{aligned}
p_{2t+2,0}
&=q^2p_{2t,0}+q^{4t+4} \\
&=q^{2t+2}\bigl([t+1]_{q^2}+q^{2t+2}\bigr)
=q^{2t+2}[t+2]_{q^2},
\end{aligned}
\]
and
\[
\begin{aligned}
p_{2t+3,0}
&=q^2p_{2t+1,0}+q^{4t+6} \\
&=q^{2t+6}\bigl([t]_{q^2}+q^{2t}\bigr)
=q^{2t+6}[t+1]_{q^2}.
\end{aligned}
\]
\end{proof}

\begin{thm}\label{thm:SO4-closed-form}
For $k\geq\ell\geq0$,
\[
p_{k,\ell}=q^\ell p_{k-\ell,0}+q^{k+1}[k-\ell+2]_q[\ell]_q+q^{k+4}(1+q^{k-\ell})\begin{bmatrix}\ell\\2\end{bmatrix}_q.
\]
\end{thm}

\begin{proof}
Iterating Theorem~\ref{thm:intro-2-shift} from $(k-\ell,0)$ gives
\[
p_{k,\ell}=q^\ell p_{k-\ell,0}+q^{k+1}\sum_{j=0}^{\ell-1}q^j\bigl([k-\ell+j+2]_q+q^2[j]_q\bigr).
\]
Using
\[
[k-\ell+j+2]_q=[k-\ell+2]_q+q^{k-\ell+2}[j]_q,\qquad
\sum_{j=0}^{\ell-1}q^j[j]_q=q\begin{bmatrix}\ell\\2\end{bmatrix}_q,
\]
we obtain
\[
p_{k,\ell}=q^\ell p_{k-\ell,0}+q^{k+1}[k-\ell+2]_q[\ell]_q
+q^{k+4}(1+q^{k-\ell})\begin{bmatrix}\ell\\2\end{bmatrix}_q.
\]
\end{proof}

Parity enters only through the boundary term in Lemma~\ref{lemma:boundary-f-values}. Theorem~\ref{thm:SO4-closed-form} thus gives a positive expression for the multiplicities computed in \cite[Theorem~7]{WallachWillenbring2000}.

\paragraph{AI disclosure.}
The first author used OpenAI's ChatGPT during manuscript preparation for editorial assistance, including improving clarity and concision, \LaTeX{} formatting, and checking proof exposition and algebraic manipulations. All mathematical arguments, references, and conclusions were reviewed by the first author, who takes full responsibility for the content of the paper.


\begin{thebibliography}{9}

\bibitem{Frohmader2026}
Andrew Frohmader, \emph{Graded multiplicities in the Kostant-Rallis setting},
Journal of Algebra \textbf{713} (2027), 592--663.

\bibitem{GoodmanWallach2009}
Roe Goodman and Nolan R. Wallach, \emph{Symmetry, Representations, and Invariants}, Graduate Texts in Mathematics, Springer, 2009.

\bibitem{KostantRallis1971}
Bertram Kostant and Stephen Rallis, \emph{Orbits and representations associated with symmetric spaces}, American Journal of Mathematics \textbf{93} (1971), no.~3, 753--809.

\bibitem{WallachWillenbring2000}
Nolan R. Wallach and Jeb Willenbring, \emph{On some $q$-analogs of a theorem of Kostant--Rallis}, Canadian Journal of Mathematics \textbf{52} (2000), no.~2, 438--448.

\bibitem{Willenbring2000}
Jeb Willenbring, \emph{Stability properties for $q$-multiplicities and branching formulas for representations of the classical groups}, Ph.D. dissertation, University of California, San Diego, 2000.

\end{thebibliography}
\end{document}